\documentclass[5p,times,authoryear]{elsarticle}

\usepackage{setspace}
\usepackage{graphicx}

\graphicspath{{./eps2pdf/}}

\usepackage{lipsum} 
\usepackage{amsmath,amsthm,amsfonts, amssymb}
\usepackage{bbm} 
\usepackage{subcaption}

\usepackage{algorithm,algorithmic}%

\usepackage{xr}
\usepackage{tikz,times,adjustbox}
\usetikzlibrary{mindmap,trees,backgrounds}
\usepackage{color}
\definecolor{red1}{rgb}{0.502,0,0}
\definecolor{blue1}{rgb}{0.098,0.098,0.498}
\definecolor{brown1}{rgb}{0.545,0.271,0.075}
\definecolor{green1}{rgb}{0,0.392,0}

\definecolor{red2}{rgb}{0.698,0.1333,0.1333}
\definecolor{blue2}{rgb}{0.275,0.51,0.706} 
\definecolor{brown2}{rgb}{0.824,0.412,0.118} 
\definecolor{green2}{rgb}{0.180,0.545,0.341}

\definecolor{red3}{rgb}{0.8039,0.3608,0.3608}
\definecolor{blue3}{rgb}{0.1176,0.565,1} 
\definecolor{brown3}{rgb}{1.0000,0.5490,0} 
\definecolor{green3}{rgb}{0.2353,0.7020,0.4431}

\theoremstyle{plain}

\theoremstyle{definition}

\theoremstyle{remark}

\DeclareMathOperator*{\esssup}{ess\,sup}

\newcommand{\floor}[1]{\left\lfloor #1 \right\rfloor}
\newcommand{\ceil}[1]{\left\lceil #1 \right\rceil}

\usepackage{multirow,array,float}

\newcolumntype{P}[1]{>{\centering\arraybackslash}p{#1}}
\newcolumntype{L}[1]{>{\flushleft\arraybackslash}p{#1}}

\journal{Annual Reviews in Control}

\begin{document}
\sloppy
{\onecolumn \large 
\section*{About} 
\begin{itemize}
	\item This is Part 2 of a two-part review paper titled ``Data-driven Decision Making in Power Systems with Probabilistic Guarantees: Theory and Applications of Chance-constrained Optimization'' by Xinbo Geng and Le Xie, Annual Reviews in Control (under review).
	\item Part 1 ``Data-driven Decision Making with Probabilistic Guarantees (Part I): A Schematic Overview of Chance-constrained Optimization'' is available at arXiv:1903.10621.
	\item Part 2 ``Data-driven Decision Making in with Probabilistic Guarantees (Part II): Applications of Chance-constrained Optimization in Power Systems'' is available at arXiv.
	\item The Matlab Toolbox \emph{ConvertChanceConstraint} (CCC) is available at \url{https://github.com/xb00dx/ConvertChanceConstraint-ccc}.
\end{itemize}
Please let us know if we missed any critical references or you found any mistakes in the manuscript.
\section*{Recent Updates} 
\begin{description}
	\item [04/2019] Part II uploaded to Arxiv.
	\item [04/2019] More CVaR-based (Convex Approximation) results are added in Part 1.
	\item [02/2019] Toolbox published at \url{https://github.com/xb00dx/ConvertChanceConstraint-ccc}. We are still working on the toolbox website and documents.
\end{description}
}

\begin{frontmatter}

\title{Data-driven Decision Making with Probabilistic Guarantees (Part II):\\
Applications of Chance-constrained Optimization in Power Systems}

\author{Xinbo Geng, Le Xie}
\address{Texas A\&M University, College Station, TX, USA.}




\begin{abstract}
Uncertainties from deepening penetration of renewable energy resources have posed critical challenges to the secure and reliable operations of future electric grids.
Among various approaches for decision making in uncertain environments, this paper focuses on chance-constrained optimization, which provides explicit probabilistic guarantees on the feasibility of optimal solutions.
Although quite a few methods have been proposed to solve chance-constrained optimization problems, there is a lack of comprehensive review and comparative analysis of the proposed methods. Part I of this two-part paper reviews three categories of existing methods to chance-constrained optimization: (1) scenario approach; (2) sample average approximation; and (3) robust optimization based methods. Data-driven methods, which are not constrained by any particular distributions of the underlying uncertainties, are of particular interest. Part II of this two-part paper provides a literature review on the applications of chance-constrained optimization in power systems. Part II also provides a critical comparison of existing methods based on numerical simulations, which are conducted on standard power system test cases.
\end{abstract}

\begin{keyword}
data-driven, power system, chance constraint, probabilistic constraint, stochastic programming, robust optimization, chance-constrained optimization.
\end{keyword}

\end{frontmatter}


\section{Introduction}
Real-time decision making in the presence of uncertainties is a classical problem that arises in many contexts. In the context of electric energy systems, a pivotal challenge is how to operate a power grid with an increasing amount of supply and demand uncertainties. The unique characteristics of such operational problem include (1) the underlying distribution of uncertainties is largely unknown (e.g. the forecast error of demand response); (2) decisions have to be made in a timely manner (e.g. a dispatch order needs to be given by 5 minutes prior to the real-time); and (3) there is a strong desire to know the risk that the system is exposed to after a decision is made (e.g. the risk of violating transmission constraints after the real-time market clears). In response to these challenges, a class of optimization problems named ``chance-constrained optimization'' has received increasing attention in both operations research and practical engineering communities. 

The objective of this article is to provide a comprehensive and up-to-date literature review on the engineering implications of chance-constrained optimization in the context of electric power systems. 

\subsection{Contributions of This Paper} 
\label{sub:contributions_of_this_paper}
The main contributions of this paper are threefold:
\begin{enumerate}
	\item We provide a detailed tutorial on existing algorithms to solve chance-constrained programs and a survey of major theoretical results. To the best of our knowledge, there is no such review available in the literature;
	\item We provide a comprehensive review on the applications of chance-constrained optimization in power systems, with focus on various interpretations of chance constraints in the context of power engineering. 
	\item We implement all the reviewed methods and develop an open-source Matlab toolbox (ConvertChanceConstraint), which is available on Github \footnote{github.com/xb00dx/ConvertChanceConstraint-ccc}. We also provide a critical comparison of existing methods based numerical simulations on IEEE standard test systems.
\end{enumerate}

\subsection{Organization of This Paper} 
\label{sub:organization_of_this_paper}
The remainder of this paper is organized as follows. Section \ref{sec:applications_in_power_systems} provides a comprehensive review on applications of CCO in power systems. The structure and usage of the Toolbox \emph{ConvertChanceConstraint} is in Section \ref{sec:numerical_simulations}.
Section \ref{sec:numerical_simulations} also conducts numerical simulations and compares existing approaches to solving CCO problems. Concluding remarks are in Section \ref{sec:concluding_remarks}.

\subsection{Notations} 
\label{sub:notations}
The notations in this paper are standard. All vectors and matrices are in the real field $\mathbf{R}$. Sets are in calligraphy fonts, e.g. $\mathcal{S}$. The upper and lower bounds of a variable $x$ are denoted by $\overline{x}$ and $\underline{x}$. The estimation of a random variable $\epsilon$ is $\hat{\epsilon}$. We use $\mathbf{1}_n$ to denote an all-one vector in $\mathbf{R}^n$, the subscript $n$ is sometimes omitted for simplicity. The absolute value of vector $x$ is $|x|$, and the cardinality of a set $\mathcal{S}$ is $|\mathcal{S}|$. Function $[a]_+$ returns the positive part of variable $a$. The indicator function $\mathbbm{1}_{x>0}$ is one if $x > 0$. The floor function $\floor{a}$ returns the largest integer less than or equal to the real number $a$. 
The ceiling function $\ceil{a}$ returns the smallest integer greater than or equal to $a$.
$\mathbb{E}[\xi]$ is the expectation of a random vector $\xi$, $\mathbb{V}(x)$ denotes the violation probability of a candidate solution $x$, and $\mathbb{P}_\xi(\cdot)$ is the probability taken with respect to $\xi$.
The transpose of a vector $a$ is $a^\intercal$. Infimum, supremum and essential supremum are denoted by $\inf$, $\sup$ and $\esssup$. The element-wise multiplication of the same-size vectors $a$ and $b$ is denoted by $a \circ b$.

\section{Applications in Power Systems} 
\label{sec:applications_in_power_systems}
A pivotal task in modern power system operation is to maintain the real-time balance of supply and demand while ensuring the system is low-cost and reliable. This pivotal task, however, faces critical challenges in the presence of rapid growth of renewable energy resources. Chance-constrained optimization, which explicitly models the risk that the system is exposed to, is a suitable conceptual framework to ensure the security and reliability of a power system under uncertainties.

\begin{figure*}[htbp]
	\centering
	\includegraphics[width=\linewidth]{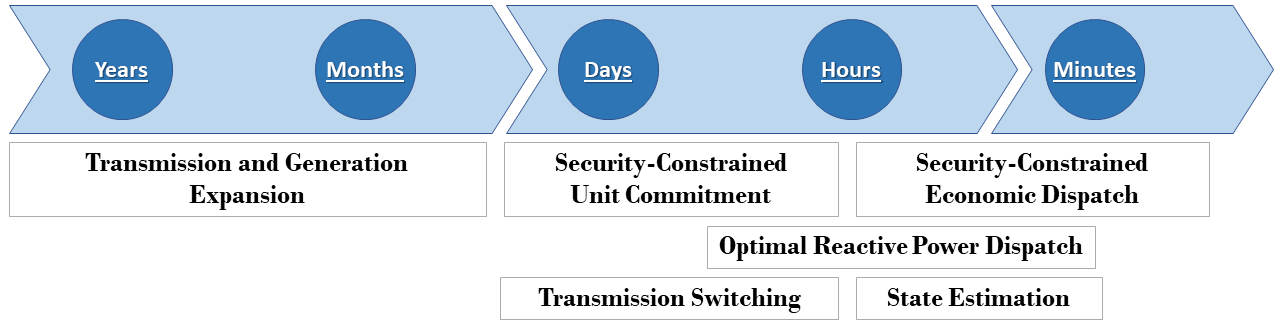}
	\caption{Representative feed-forward decisions made in power system planning and operation}
	\label{fig:feed_forward_decision_makin_framework}
\end{figure*}

There is a large body of literature adopting CCO for power system applications. Figure \ref{fig:feed_forward_decision_makin_framework} presents some existing applications of CCO in power systems. In the following sections, we introduce three important applications of CCO in power systems: security-constrained economic dispatch (SCED) (Section \ref{sub:security_constrained_economic_dispatch}), security-constrained unit commitment (SCUC) (Section \ref{sub:security_constrained_unit_commitment}) and generation and transmission expansion (Section \ref{sub:generation_and_transmission_expansion}).

Figure \ref{fig:feed_forward_decision_makin_framework} also presents a feed-forward decision making framework for power system operations. The feed-forward framework partitions the overall decision making process into several time segments. The longer-term decisions (e.g. generation expansion) are fed into shorter-term decision making processes (e.g. unit commitment). The shorter-term decisions (e.g. generation commitment from SCUC) have direct impacts on real-time operations (e.g. dispatch results in SCED). It is worth noting that as time draws closer to the actual physical operation, information gets much sharper \citep{xie_wind_2011}. 

%


\begin{table*}[tb]
	\caption{Power System Applications of Chance-constrained Optimization}
	\label{tab:cco_applications_in_power_systems}
	\centering

	\begin{tabular}{P{0.04\linewidth}|L{0.1\linewidth}L{0.15\linewidth}L{0.15\linewidth}L{0.3\linewidth}L{0.1\linewidth}}
	\hline

	\hline
	 & \textbf{Methods} & \textbf{Expansion} & \textbf{SCUC} & \textbf{SCED} & \textbf{Other Applications} \\
	\hline
	\multirow{1}{*}{\rotatebox[origin=c]{90}{\textbf{Deterministic}}\rotatebox[origin=c]{90}{\textbf{Equivalent}}} & Gaussian
	& \citep{sanghvi_strategic_1982}\citep{lopez_generation_2007}\citep{mazadi_modified_2009}\cite{manickavasagam_sensitivity-based_2015}
	& \citep{ding_studies_2010,pozo_chance-constrained_2013,wu_chance-constrained_2014}
	& \citep{bent_synchronization-aware_2013,bienstock_robust_2013,jabr_adjustable_2013,roald_analytical_2013,bienstock_chance-constrained_2014,li_analytical_2015,roald_optimal_2015,zhang_data-driven_2015,lubin_robust_2016,doostizadeh_energy_2016,roald_corrective_2017,roald_optimized_2017,wang_chance-constrained_2017,vrakopoulou_chance_2019,li_chance_2019}
	& \citep{lopez_convex_2015,franco_robust_2016}\\
	\hline
	\multirow{2}{*}{\rotatebox[origin=c]{90}{\textbf{Scenario}}\rotatebox[origin=c]{90}{\textbf{Approach}}}  & a-priori 
	& -
	& \citep{geng_security-constrained_2019}
	& \citep{bucher_probabilistic_2013,vrakopoulou_probabilistic_2013,vrakopoulou_probabilistic_2013-3,roald_risk-constrained_2014,roald_risk-based_2015,zhang_data-driven_2015,ming_scenario-based_2017,modarresi_scenario-based_2018,geng_data-driven_2019-1}
	& \citep{yang_joint_2014}\\
	\cline{2-6}
	& a-posteriori
	& -
	& \citep{margellos_stochastic_2013,geng_security-constrained_2019,hreinsson_stochastic_2015} 
	& \citep{modarresi_scenario-based_2018,geng_data-driven_2019-1} & - \\
	\hline
	\multirow{1}{*}{\rotatebox[origin=c]{90}{\textbf{Sample}}\rotatebox[origin=c]{90}{\textbf{Average}}\rotatebox[origin=c]{90}{\textbf{Approximation}}}  & - 
	& \citep{zhang_chance-constrained_2017}
	& \citep{wang_chance-constrained_2012,wang_price-based_2013,zhao_expected_2014,tan_hybrid_2016,bagheri_data-driven_2017,zhang_chance-constrained_2017-1}
	& \citep{geng_data-driven_2019-1} & - \\ 
	\hline
	\multirow{2}{*}{\rotatebox[origin=c]{90}{\textbf{RO-based}}\rotatebox[origin=c]{90}{\textbf{Approach}}}
	& RLO & - & \cite{jiang_robust_2012} & \citep{geng_data-driven_2019-1} & - \\
	\cline{2-6}
	& Convex Approximation
	& - 
	& - 
	& \citep{zhang_robust_2013,summers_stochastic_2014,summers_stochastic_2015,geng_data-driven_2019-1}
	& - \\
	\hline 
	\multirow{1}{*}{\rotatebox[origin=c]{90}{\textbf{Others}}} & - & \citep{yang_chance_2005,qiu_risk-based_2016} & \citep{martinez_risk-averse_2015,wu_solution_2016} & \citep{vrakopoulou_probabilistic_2013,bienstock_chance-constrained_2014,doostizadeh_energy_2016,ke_novel_2016,muhlpfordt_solving_2017,wang_chance-constrained_2017} & -\\
	\hline

	\hline
	\end{tabular}
\end{table*}

\subsection{Security-Constrained Economic Dispatch} 
\label{sub:security_constrained_economic_dispatch}
\subsubsection{Deterministic SCED} 
\label{ssub:deterministic_sced}
Security-constrained Economic Dispatch (SCED) lies at the center of modern electricity markets and short-term power system operations. It determines the most cost-efficient output levels of generators while keeping the real-time balance between supply and demand. Different variations of the SCED problem are all based on the direct current optimal power flow (DCOPF) problem. We present a typical form of DCOPF with wind generation.
\begin{subequations}
\label{form:det_dcopf}
  \begin{align}
    \text{(det-DCOPF):}~\min_g~& c(g) \label{form:det-DCOPF-obj} \\
    \text{s.t}~& \mathbf{1}^\intercal g = \mathbf{1}^\intercal d - \mathbf{1}^\intercal \hat{w} \label{form:det-DCOPF-balance} \\ 
    & f = H_g g + H_w \hat{w} - H_d d \label{form:det-DCOPF-flow} \\ 
    & \underline{f} \le f \le \overline{f} \label{form:det-DCOPF-flow-limits} \\
    & \underline{g} \le g \le \overline{g} \label{form:det-DCOPF-gen-limits}
  \end{align}
\end{subequations}
The decision variables are generation output levels $g \in \mathbf{R}^{n_g}$. The objective of (det-DCOPF) is to minimize total generation cost $c(g)$, while ensuring total generation equates total \emph{net} demand \footnote{Wind generation is treated as negative loads.} (\ref{form:det-DCOPF-balance}). Constraints include transmission line flow limits (\ref{form:det-DCOPF-flow})-(\ref{form:det-DCOPF-flow-limits}) and generation capacity limits (\ref{form:det-DCOPF-gen-limits}). Transmission line flows $f \in \mathbf{R}^{n_l}$ are calculated using (\ref{form:det-DCOPF-flow}), in which $H$ is the power transfer distribution factor (PTDF) matrix, and $H_g \in \mathbf{R}^{n_l \times n_g}$ ($H_d \in \mathbf{R}^{n_l \times n_d}$,$H_w \in \mathbf{R}^{n_l \times n_w}$) denotes the submatrix formed by the columns of $H$ corresponding to generators (loads, wind farms). (\ref{form:det_dcopf}) utilizes the expected wind generation or wind forecast $\hat{w}$, we refer to (\ref{form:det_dcopf}) as \emph{deterministic DCOPF} (det-DCOPF) since no uncertainties are being considered. 

\subsubsection{Chance-constrained SCED} 
\label{ssub:chance_constrained_sced}
Many researchers advance (det-DCOPF) towards a chance-constrained formulation with wind uncertainties. A representative formulation is (\ref{form:cc-DCOPF}), which appears in a majority of the existing literatures, e.g. \citep{bienstock_chance-constrained_2014,vrakopoulou_probabilistic_2013}.
\begin{subequations}
\label{form:cc-DCOPF}
  \begin{align}
  & \text{(cc-DCOPF):} \nonumber \\
    \min_{g,\eta}~& c(g) \\
    \text{s.t}~ & \mathbf{1}^\intercal g = \mathbf{1}^\intercal d - \mathbf{1}^\intercal \hat{w} \label{form:cc-DCOPF-balance} \\
    & f(\hat{w},\tilde{w}) =  H_g (g - \mathbf{1}^\intercal \tilde{w} \eta )\nonumber\\
    & \hspace{40pt} - H_d d + H_w (\hat{w} + \tilde{w}) \\
    & \mathbb{P}_{\tilde{w}} \Big ( \underline{f} \le f(\hat{w},\tilde{w})  \le \overline{f}~\text{and}\nonumber\\
    & \hspace{15pt} \underline{g} \le g - \mathbf{1}^\intercal \tilde{w} \eta \le \overline{g} \Big ) \ge 1 - \epsilon \label{form:cc-DCOPF-cc} \\
    & \mathbf{1}^\intercal \eta = 1 \label{form:cc-DCOPF-sumone}\\    
    & \underline{g} \le g \le \overline{g} \\
    & - \mathbf{1} \le \eta \le \mathbf{1}
  \end{align}
\end{subequations}
Unlike (det-DCOPF) using wind forecast $\hat{w}$, chance-constrained SCED (cc-DCOPF) explicitly models wind generation as a random vector $w \in \mathbf{R}^{n_w}$. The wind generation $w = \hat{w}+ \tilde{w}$ is decomposed into two components: the \emph{deterministic} wind forecast value $\hat{w} \in \mathbf{R}^{n_w}$ and the \emph{uncertain} forecast error $\tilde{w} \in \mathbf{R}^{n_w}$.
To guarantee the real-time balance of supply and demand, (cc-DCOPF) introduces an affine control policy $\eta \in [-1,1]^{n_g}$ to proportionally allocate total wind fluctuations $\mathbf{1}^\intercal \tilde{w}$ to each generator. It is easy to verify that constraints (\ref{form:cc-DCOPF-balance}) and (\ref{form:cc-DCOPF-sumone}) imply the supply-demand balance in the presence of wind uncertainties, i.e.
\begin{equation}
	\mathbf{1}^\intercal (g - \mathbf{1}^\intercal \tilde{w} \eta ) = \mathbf{1}^\intercal d - \mathbf{1}^\intercal ( \hat{w} +\tilde{w}),
\end{equation}
The affine policy vector $\eta \in \mathbf{R}^{n_g}$ is sometimes referred as participation factor or distribution vector \citep{vrakopoulou_probabilistic_2013}.
The (joint) chance constraint (\ref{form:cc-DCOPF-cc}) constrains the transmission flow and generation within their capacities with high probability $1- \epsilon$ in the presence of wind uncertainties.

For simplicity, we only account for the major source of uncertainties (i.e. wind) in the real-time. Many references provides more complicated formulation of (cc-DCOPF), e.g. considering joint uncertainties from load and wind \citep{doostizadeh_energy_2016,muhlpfordt_solving_2017}, and contingencies of potential generator or transmission line outages \citep{roald_optimal_2015}.

There exist a few different but similar formulations of (cc-DCOPF). In general, policies of any form could help balance supply with demand under uncertainties. The affine policy in (cc-DCOPF) is the simplest choice and lead to optimization problems that are easy to solve. There are other papers applying different forms of policies, e.g. \citep{jabr_adjustable_2013} introduces a matrix form of the affine policy $\Upsilon \in \mathbf{R}^{n_g \times n_w}$, which specifies the corrective control of each generator on each wind farm.
(cc-DCOPF) is a single snapshot dispatch problem, it is straightforward to extend it to a multi-period or look-ahead dispatch problem \citep{modarresi_scenario-based_2018,vrakopoulou_probabilistic_2013}. Many papers evaluate the impacts of new elements in modern power systems, such as demand response \citep{ming_scenario-based_2017,zhang_distributionally_2017}, ambient temperatures and meteorological quantities \citep{bucher_probabilistic_2013}, and frequency control \citep{li_analytical_2015,zhang_distributionally_2017}.

\subsubsection{Solving cc-DCOPF} 
\label{ssub:solving_cc_sced}
Table \ref{tab:cco_applications_in_power_systems} summarizes various methods to solve (cc-DCOPF). The most popular one consists of two steps: (i) decomposing the joint chance constraint (\ref{form:cc-DCOPF-cc}) into individual ones  $\mathbb{P}_\xi(f_i(x,\xi) \le 0) \ge 1- \epsilon_i, i=1,2,\cdots,m$; (ii) deriving the deterministic equivalent form of each individual chance constraint by making the Gaussian assumption. More technical details of this method are in Section \ref{sub:special_cases} of \citep{geng_data-driven_2019}. This method is taken by many researchers for its simplicity and computationally tractable reformulation. Although the Gaussian assumption enjoys the law of large numbers, it is often an approximation or even doubtful assumption. For example, \citep{hodge_wind_2011} shows that the wind forecast error is better represented by Cauchy distributions instead of Gaussian ones. The first step of this method is to decompose a joint chance constraint $\mathbb{P}_\xi(f(x,\xi) \le 0) \ge 1- \epsilon$ into individual ones. As discussed in Section \ref{sub:joint_and_individual_chance_constraints} and \ref{ssub:converting_a_joint_chance_constraint_to_individual} of \citep{geng_data-driven_2019}, this step often introduces conservativeness because of the limitation of Bonferroni inequality. The level of conservativeness could be significant when the number of constraints $m$ is large, which is typically the case in power systems. 

The scenario approach is another commonly-accepted method. It provides rigorous guarantees on the quality of the solution and does not assume the distribution is Gaussian or any particular type. Most papers adopting the scenario approach apply the \emph{a-priori} guarantees (e.g. Theorem \ref{thm:exact_feasibility_scenario_approach} and \ref{thm:prior_guarantee_helly_dim} in \citep{geng_data-driven_2019}) on (cc-DCOPF) and verify the a-posteriori feasibility of solutions through Monte-Carlo simulations. One common observation is that the solution $x_N^\ast$ is often quite conservative, i.e. $\mathbb{V}(x_N^\ast) \ll \epsilon$.
One major source of conservativeness is the loose sample complexity bounds $N$ \footnote{Many papers still utilize the first sample complexity bound proved in \citep{calafiore_uncertain_2005}, which was significantly tightened in \citep{campi_exact_2008} and following works \citep{calafiore_random_2010}.}.
Since (cc-DCOPF) is convex, Theorem \ref{thm:max_num_support_scenario} in \citep{geng_data-driven_2019} states that the number of decision variables $n$ is an upper bound of the number of support scenarios $|\mathcal{S}|$ or Helly's dimension $h$. This upper bound, as pointed out in \citep{modarresi_scenario-based_2018}, is indeed very loose. \citep{modarresi_scenario-based_2018} reported only $\sim 5$ support scenarios for a chance-constrained look-ahead SCED problem with thousands of decision variables. By exploiting the structural features of (cc-DCOPF), the sample complexity bound $N$ can be significantly improved. Unfortunately, only \citep{modarresi_scenario-based_2018} and \citep{ming_scenario-based_2017} followed this path to reduce conservativeness. 

There are also many papers utilizing the robust optimization related methods to solve (cc-DCOPF). \citep{jiang_robust_2012} constructs uncertainty sets with the help of probabilistic guarantees in \citep{bertsimas_price_2004}. References \citep{summers_stochastic_2014,summers_stochastic_2015} incorporate the convex approximation framework and compare different choices of generating functions $\phi(z)$ on (cc-DCOPF). Although there are no explicit forms of chance constraints in \citep{zhang_robust_2013}, the CVaR-oriented approach therein can be interpreted as solving cc-DCOPF using convex approximation with the choice of Markov bound. 

Most papers in Table \ref{tab:cco_applications_in_power_systems} aim at finding suboptimal solutions to (cc-DCOPF). However, it is somewhat surprising to note that none of them estimates how suboptimal the solution is via approaches like Proposition \ref{prop:lower_bound_cco_framework} or \ref{prop:lower_bound_cco_saa} in \citep{geng_data-driven_2019}. Almost all the papers evaluate the a-posteriori feasibility by Monte-Carlo simulations with a huge sample size. Methods like Proposition \ref{prop:check_feasibility_scenario_approach}  in \citep{geng_data-driven_2019} would be more attractive when data is limited, which is closer to the reality.

\subsection{Security-Constrained Unit Commitment} 
\label{sub:security_constrained_unit_commitment}
\subsubsection{Deterministic SCUC} 
\label{ssub:deterministic_scuc}
Security-Constrained Unit Commitment (SCUC) is one of the most important procedures in power system day-ahead or intra-day operations.
\begin{subequations}
\label{opt:det-SCUC}
\begin{align}
& \text{(det-SCUC):} \nonumber \\
  \min_{z,u,v,g,s}~ & \sum_{t=1}^{n_t} c_n^\intercal z^t + c_u^\intercal u^t + c_v^\intercal v^t + c_g^\intercal g^{t,0} + c_s^\intercal s^t \label{opt:det-SCUC-obj} \\
\text{s.t.}~& \mathbf{1}^\intercal g^{t,k} \ge \mathbf{1}^\intercal \hat{d}^t - \mathbf{1}^\intercal \hat{w}^t \label{opt:det-SCUC-balance}\\
& \underline{f} \le H_g^{t,k} g^{t,k} - H_d^{t,k} \hat{d}^t + H_w^{t,k} \hat{w}^t \le \overline{f} \label{opt:det-SCUC-line} \\
& \underline{r} \le g^{t,k} - g^{t-1,k} \le \overline{r}  \label{opt:det-SCUC-ramp} \\
& a^k \circ (g^{t,0} - s^t) \le g^{t,k} \le a^k \circ (g^{t,0} + s^t) \label{opt:det-SCUC-contingency} \\
& \hspace{114pt} k \in [0,n_k], t \in [1,n_t] \nonumber \\
& \underline{g} \circ  z^t \le g^{t,0} \le \overline{g} \circ  z^t \label{opt:det-SCUC-gencap}\\
& \underline{s} \circ z^t \le s^t \le \overline{s} \circ  z^t \label{opt:det-SCUC-reservecap} \\
& \underline{g} \circ  z^t \le g^{t,0} - s^t \le g^{t,0} + s^t \le \overline{g} \circ  z^t \label{opt:det-SCUC-resgencap}\\
& z^{t-1} - z^t + u^t \ge 0 \label{opt:det-SCUC-startup} \\
& z^t - z^{t-1} + v^t \ge 0 \label{opt:det-SCUC-shutdown} \\
& \hspace{162pt} t \in [1,n_t] \nonumber \\
& z_i^{t} - z_i^{t-1} \le z_i^\iota,~\iota \in [t+1,\min\{t+\underline{u}_i-1,n_t\}] \label{opt:det-SCUC-minon-time} \\
& z_i^{t-1} - z_i^{t} \le 1- z_i^\iota, ~ \iota \in [t+1,\min\{t+\underline{v}_i-1,n_t\}] \label{opt:det-SCUC-minoff-time}\\
& \hspace{110pt} i \in [1,n_g],~t \in [2,n_t]  \nonumber
\end{align}
\end{subequations}
Deterministic SCUC (det-SCUC) seeks the optimal commitment and generation schedule of $n_g$ generators for the upcoming $n_t$ snapshot while ensuring system security in $n_k$ contingencies. Decision variables include commitment and startup/shutdown decisions $(z^t,u^t, v^t)$, as well as generation and reserve schedules $(g^{t,k}, s^t)$. 
The objective of (\ref{opt:det-SCUC}) is to minimize total operation costs, which include 
no-load costs $c_n^\intercal z^t$, startup costs $c_u^\intercal u^t$, shutdown costs $c_v^\intercal v^t$, generation costs $c_g^\intercal g^{t,0}$ and reserve costs $c_s^\intercal s^t$. Constraint (\ref{opt:det-SCUC-balance}) assures there is enough supply to meet \emph{net} demand. Constraints (\ref{opt:det-SCUC-line}), (\ref{opt:det-SCUC-ramp}) and (\ref{opt:det-SCUC-reservecap}) are about transmission capacity, generation ramping capability and reserve limit in contingency scenario $k$ at time $t$. In contingency scenarios, the adjusted output $g_i^{t,k}$ of generator $i$ is bounded by its reserve $s_i^t$. Vector $a^k \in \{0,1\}^{n_g}$ represents the availability of generators in contingency $k$. When $a_i^k = 0$, generator $i$ is available in contingency $k$, thus has zero generation output.
Generation and reserve capacity constraints are in (\ref{opt:det-SCUC-gencap}) and (\ref{opt:det-SCUC-reservecap}). Constraints (\ref{opt:det-SCUC-gencap})-(\ref{opt:det-SCUC-resgencap}) also ensure the consistency of generation with commitment decisions. (\ref{opt:det-SCUC-startup})-(\ref{opt:det-SCUC-shutdown}) are the logistic constraints about commitment status, startup and shutdown decisions. Minimum on/off time constraints for all generators are presented in (\ref{opt:det-SCUC-minon-time})-(\ref{opt:det-SCUC-minoff-time}).
\subsubsection{Chance-constrained SCUC} 
\label{ssub:chance_constrained_scuc}
Many researchers proposed various advanced formulations of SCUC to deal with uncertainties, e.g. using robust optimization \citep{bertsimas_adaptive_2013} and stochastic programming \citep{takriti_stochastic_1996}. A good overview of SCUC formulations with uncertainties is in \citep{zheng_stochastic_2015}. In this paper, we formulate the chance-constrained SCUC problem.
Unlike the case of SCED, there is no unified formulation of chance-constrained SCUC.
We present one simplified formulation in (\ref{opt:cc-SCUC}). Alternative formulations of chance-constrained SCUC can be found in \citep{jiang_robust_2012,wu_stochastic_2007,zheng_stochastic_2015}.
\begin{subequations}
\label{opt:cc-SCUC}
\begin{align}
& \text{(cc-SCUC):} \nonumber \\
  \min_{z,u,v,g,s}~ & \sum_{t=1}^{n_t} c_n^\intercal z^t + c_u^\intercal u^t + c_v^\intercal v^t + c_g^\intercal g^{t,0} + c_s^\intercal s^t \label{opt:cc-SCUC-obj} \\
\text{s.t.}~& (\ref{opt:det-SCUC-balance}), (\ref{opt:det-SCUC-line}), (\ref{opt:det-SCUC-ramp}), (\ref{opt:det-SCUC-contingency}),~ k \in [0,n_k], t \in [1,n_t] \nonumber \\
&  (\ref{opt:det-SCUC-gencap}), (\ref{opt:det-SCUC-reservecap}), (\ref{opt:det-SCUC-resgencap})),(\ref{opt:det-SCUC-startup}), (\ref{opt:det-SCUC-shutdown}),~t \in [1,n_t] \nonumber \\
& (\ref{opt:det-SCUC-minon-time}),(\ref{opt:det-SCUC-minoff-time}),~i \in [1,n_g],~t \in [2,n_t] \nonumber \\
& \mathbb{P}\Big( \mathbf{1}^\intercal g^{t,k} \ge \mathbf{1}^\intercal (\hat{d}^t+\tilde{d}^t) - \mathbf{1}^\intercal (\hat{w}^t+\tilde{w}^t), \label{opt:cc-SCUC-balance} \\
& \quad \underline{f} \le H_g^{t,k} g^{t,k} - H_d^{t,k} (\hat{d}^t+\tilde{d}^t) \nonumber \\ 
& \hspace{55pt}+ H_w^{t,k} (\hat{w}^t+\tilde{w}^t) \le \overline{f}, \label{opt:cc-SCUC-line} \\
& \hspace{28pt} k \in [0,n_k], t \in [1,n_t] \Big) \ge 1 - \epsilon \label{opt:cc-SCUC-chance}
\end{align}
\end{subequations}
The formulation of (cc-SCUC) is almost identical to (det-SCUC) except the chance constraint (\ref{opt:cc-SCUC-balance})-(\ref{opt:cc-SCUC-chance}). In (cc-SCUC), wind generation $w \in \mathbf{R}^{n_w}$ is modeled as a random vector consisting of a deterministic predicted component $\hat{w} \in \mathbf{R}^{n_w}$ and a stochastic error component $\tilde{w} \in \mathbf{R}^{n_w}$. The chance constraint (\ref{opt:cc-SCUC-balance})-(\ref{opt:cc-SCUC-chance}) ensures enough supply to meet demand and line flows within limits under uncertainties with probability at least $1- \epsilon$ for any contingency scenario $k$ at any time $t$.

The \emph{joint} chance constraint (\ref{opt:cc-SCUC-balance})-(\ref{opt:cc-SCUC-chance}) is sometimes written as two (joint) chance constraints:
\begin{subequations}
\begin{align}
	& \mathbb{P}\Big( \mathbf{1}^\intercal g^{t,k} \ge \mathbf{1}^\intercal (\hat{d}^t+\tilde{d}^t) - \mathbf{1}^\intercal (\hat{w}^t+\tilde{w}^t),\nonumber \\
	& \hspace{60pt} k \in [0, n_k],~t\in[1,n_t] \Big) \ge 1 - \epsilon^{\text{LOLP}} \label{opt:cc-SCUC-LOLP} \\
	& \mathbb{P}\Big( \underline{f} \le H_g^{t,k} g^{t,k} - H_d^{t,k} (\hat{d}^t+\tilde{d}^t) + H_w^{t,k} (\hat{w}^t+\tilde{w}^t) \le \overline{f}, \nonumber \\
	& \hspace{60pt} k \in [0, n_k],~t\in[1,n_t] \Big) \ge 1 - \epsilon^{\text{TLOP}} \label{opt:cc-SCUC-TLOP}
\end{align}
\end{subequations}
An important metric to evaluate power system reliability is through the \emph{loss of load probability} (LOLP), which is defined as the probability that the total demand is not met by the total generation \citep{allan_reliability_2013,qiu_risk-based_2016}. It can be seen that (\ref{opt:cc-SCUC-LOLP}) is essentially ensuring the value of LOLP will not exceed a desired level $\epsilon^\text{LOLP}$. Similarly, we could define the concept \emph{transmission line overload probability} (TLOP) \citep{wu_chance-constrained_2014}. Then (\ref{opt:cc-SCUC-TLOP}) is the same as $\text{TLOP} \le \epsilon^\text{TLOP}$. 

Some papers (e.g. \citep{wu_chance-constrained_2014}) further break down the joint chance constraint (\ref{opt:cc-SCUC-LOLP})-(\ref{opt:cc-SCUC-TLOP}) into individual chance constraints (\ref{opt:cc-SCUC-LOLP-individual})-(\ref{opt:cc-SCUC-TLOP-individual}), which can be interpreted as constraints on LOLP or TLOP for each time period $t$.
\begin{subequations}
\begin{align}
	& \mathbb{P}\Big( \mathbf{1}^\intercal g^{t,k} \ge \mathbf{1}^\intercal (\hat{d}^t+\tilde{d}^t) - \mathbf{1}^\intercal (\hat{w}^t+\tilde{w}^t) \Big) \ge 1 - \epsilon_{t,k}^{\text{LOLP}},  \nonumber \\
	& \hspace{120pt} k \in [0, n_k],~t\in[1,n_t]. \label{opt:cc-SCUC-LOLP-individual} \\
	& \mathbb{P}\Big( \underline{f} \le H_g^{t,k} g^{t,k} - H_d^{t,k} (\hat{d}^t+\tilde{d}^t) + H_w^{t,k} (\hat{w}^t+\tilde{w}^t) \le \overline{f} \Big) \ge 1 - \epsilon_{t,k}^{\text{TLOP}},  \nonumber \\
	& \hspace{120pt} k \in [0, n_k],~t\in[1,n_t].  \label{opt:cc-SCUC-TLOP-individual} 
\end{align}
\end{subequations}
Another interesting set of chance constraints in cc-SCUC guarantees the utilization ratio of wind generation greater than a desired threshold with high probability $1- \epsilon$ \citep{wang_chance-constrained_2012,wang_price-based_2013,zhao_expected_2014}. Different variations of the chance constraint on wind utilization ratios can be found in \citep{wang_chance-constrained_2012}.

\subsubsection{Solving Chance-constrained SCUC} 
\label{ssub:solving_chance_constrained_scuc}
As mentioned in Section \ref{ssub:chance_constrained_scuc}, there is no uniform formulation of chance-constrained SCUC. Many references in Table \ref{tab:cco_applications_in_power_systems} concentrate on exploring alternative formulations of cc-SCUC. Therefore theoretical guarantees or quality solution is not a major concern. 

Among all the reviewed methods, sample average approximation is commonly used when solving chance-constrained SCUC \citep{wang_chance-constrained_2012,wang_price-based_2013,zhao_expected_2014,tan_hybrid_2016,bagheri_data-driven_2017,zhang_chance-constrained_2017-1}. Section \ref{sec:sample_average_approximation} of \citep{geng_data-driven_2019} shows that SAA reformulates (CCO) to a mixed integer program, which is difficult to solve in general. Many references apply various techniques from integer programming to speed up the computation, e.g. \citep{zhao_expected_2014,jiang_cutting_2016}.

Section \ref{sub:structural_properties_of_the_scenario_problem} of \citep{geng_data-driven_2019} shows that there is no upper bound on the number of support scenarios for non-convex problems in general. Thus, a majority of results of the scenario approach cannot be directly applied on cc-SCUC. Recently, \citep{campi_general_2018} extends the a-posteriori guarantees of the scenario approach towards non-convex problems. 
\citep{geng_security-constrained_2019} adopts the approach in \citep{campi_general_2018} and shows the possibility to apply the theoretical results of the scenario approach on (cc-SCUC). It is worth mentioning that some theoretical results in robust optimization still apply in spite of the non-convexity of SCUC from integer variables $(z^t,u^t,v^t)$, e.g. \citep{bertsimas_data-driven_2018}. This could be an interesting direction to explore.

\subsection{Generation and Transmission Expansion} 
\label{sub:generation_and_transmission_expansion}
Generation and transmission expansion (the expansion problem in short) is a critical component in \emph{long-term} power system planning exercises. The expansion problem answers the following critical questions: (i) when to invest on new elements such as transmission lines and generators in the system; (ii) what types of new elements are necessary; and (iii) how much capacity is needed and where the best locations would be for those new elements. A typical objective of the expansion problem is to minimize (i) total cost of investment in new generators and transmission line; (ii) environmental impacts; and (iii) cost of generation.
Constraints of the expansion problem often include total or individual costs within budget, capacity constraint, reliability requirement, supply-demand balance, power flow equations, and operation requirements such as generation or transmission limits.

The expansion problem typically needs to deal with uncertainties from demand, generation and transmission outages, and renewables. Chance constraints often appear as requirements on reliability metrics such as LOLP (\ref{opt:cc-SCUC-LOLP}) and TLOP (\ref{opt:cc-SCUC-TLOP}).

Among all the papers incorporating chance constraints in the expansion problem, a majority of them assume the underlying distribution is Gaussian and derive the second order cone equivalent form (Section \ref{sub:special_cases}, \citep{geng_data-driven_2019}), e.g. \citep{sanghvi_strategic_1982,lopez_generation_2007,mazadi_modified_2009,manickavasagam_sensitivity-based_2015}. A few papers design its own simulation-based iterative algorithms because of complicated problem formulations, e.g. \citep{yang_chance_2005,qiu_risk-based_2016}. Although Monte-Carlo simulation is typically performed to evaluate the actual feasibility, there is no rigorous guarantees on these results.

Similar to the chance constrained DCOPF problem, deriving deterministic equivalent forms is the most popular choice. Considering the expansion problem is usually ultra-large-scale and involves lots of integer variables, the simplicity of deterministic equivalent form becomes particularly attractive. Additional pros and cons of this approach are analyzed in Section \ref{ssub:solving_cc_sced}. 

Similar to chance-constrained SCUC, the expansion problem includes many integer variables and is non-convex in nature. As discussed in Section \ref{ssub:solving_chance_constrained_scuc}, the scenario approach and sample average approximation can still be applied on the expansion problem. Because of the size of the expansion problem, the required sample complexity could be astronomic, which lead to major computational issues. Although the scenario approach and sample average approximation could provide better theoretical guarantees, it is essential to overcome the major obstacles in computation to apply some better methods on the expansion problem.






\section{Numerical Simulations} 
\label{sec:numerical_simulations}

\subsection{Simulation Settings} 
\label{sub:simulation_settings}
Chance-constrained DCOPF (\ref{form:cc-DCOPF}) serves as a benchmark problem for a critical comparison of solutions to (CCO). We provide numerical solutions of cc-DCOPF on two test systems: a 3-bus system and the IEEE 24-bus RTS test system.

The 3-bus system is a modified version of the 3-bus system in \citep{lesieutre_examining_2011}. The major difference is the removal of the load at bus 2 and the synchronous condensor at bus 3 in order to visualize the feasible region and the space of uncertainties. The original 3-bus system ``\emph{case3sc.m}'' is available in the Matpower toolbox \cite{zimmerman_matpower:_2011}. The modified system in this paper can be found in the examples of CCC \footnote{github.com/xb00dx/ConvertChanceConstraint-ccc/tree/master/examples}. For simplicity, we only consider uncertainties of loads, which is modeled as Gaussian variables with 5\% standard variation. 

The 24-bus system in this paper is a modified version of the IEEE 24-bus RTS benchmark system \citep{grigg_ieee_1999}. The transmission line capacities are set to be 60\% of the original capacities.
We conduct two sets of simulations on the 24-bus system with different distributions of uncertainties. 
The first one is similar with the 3-bus case, nodal loads are modeled as independent Gaussian variables with 5\% standard deviation. The second one models the errors of nodal load forecasts as independent beta-distributed random variables, with parameters $\alpha=25.2414$ and $\beta=25.2692$ \footnote{This setting of beta distribution is from \citep{hodge_wind_2011}, and scaled from $[0,1]$ to $[-18\%,18\%]$.}.

Ten Monte-Carlo simulations are conducted on every method to examine the randomness of solutions. For the 3-bus case, each Monte-Carlo simulation uses 100 i.i.d samples to solve cc-DCOPF. 2048 points are used in each run to solve (cc-DCOPF) of the 24-bus system. The returned solutions are evaluated on an independent set of $10^4$ points. 

We use Gurobi 7.10 \citep{gurobi_optimization_gurobi_2016} to get results of scenario approach and sample average approximation. Cplex 12.8 is used to solve (CCO) with robust counterpart and convex approximation. 

\subsection{Feasible Region} 
\label{sub:case_study_feasible_region}
\begin{figure}[tb]
	\centering
	\includegraphics[width=\linewidth]{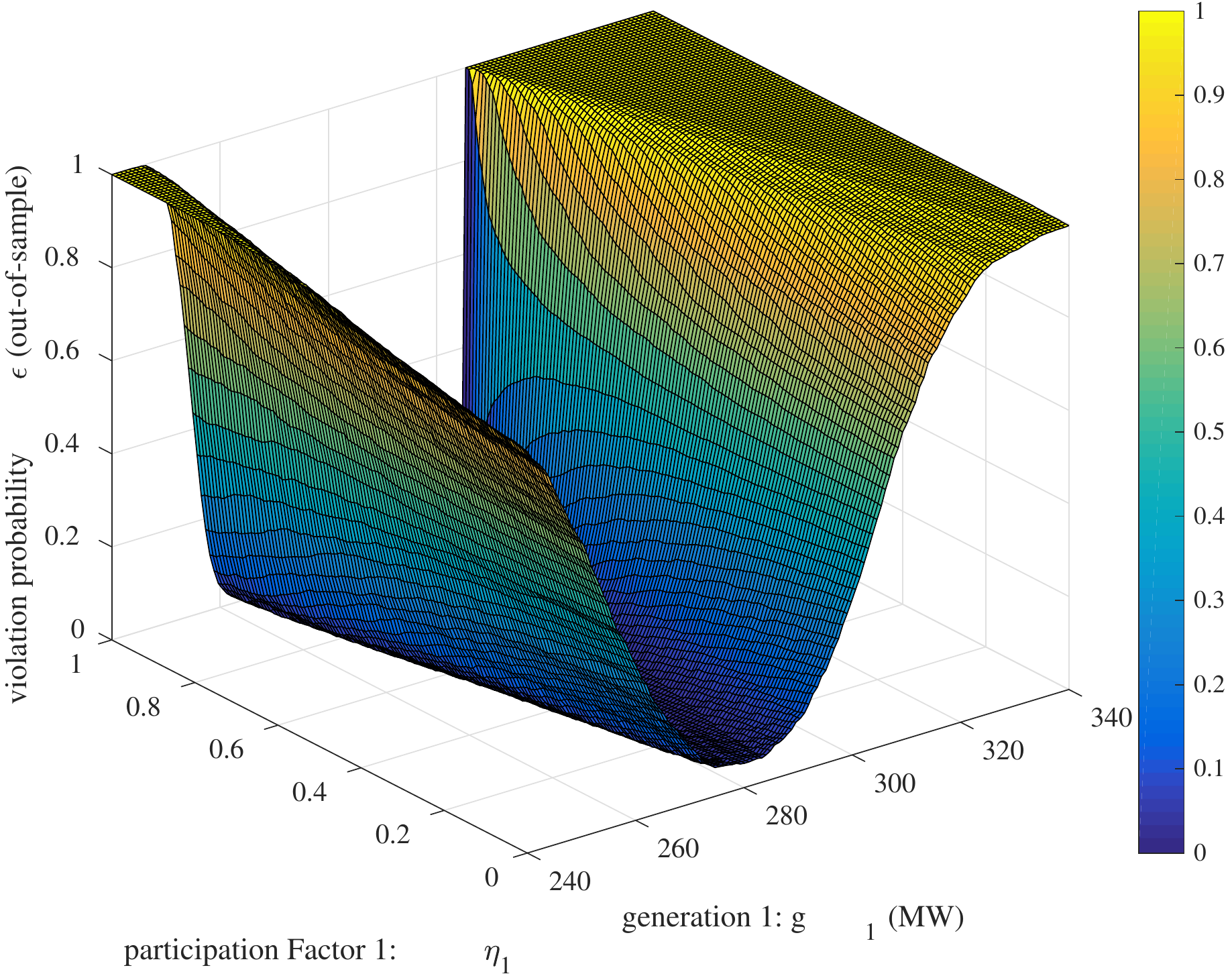}
	\caption{Visualization of the Violation Probability Function $\mathbb{V}(x)$ for cc-DCOPF of the 3-bus System}
	\label{fig:ex_case3sc-feasible-region-surf-view3}
\end{figure}

\begin{figure}[tb]
	\centering
	\includegraphics[width=\linewidth]{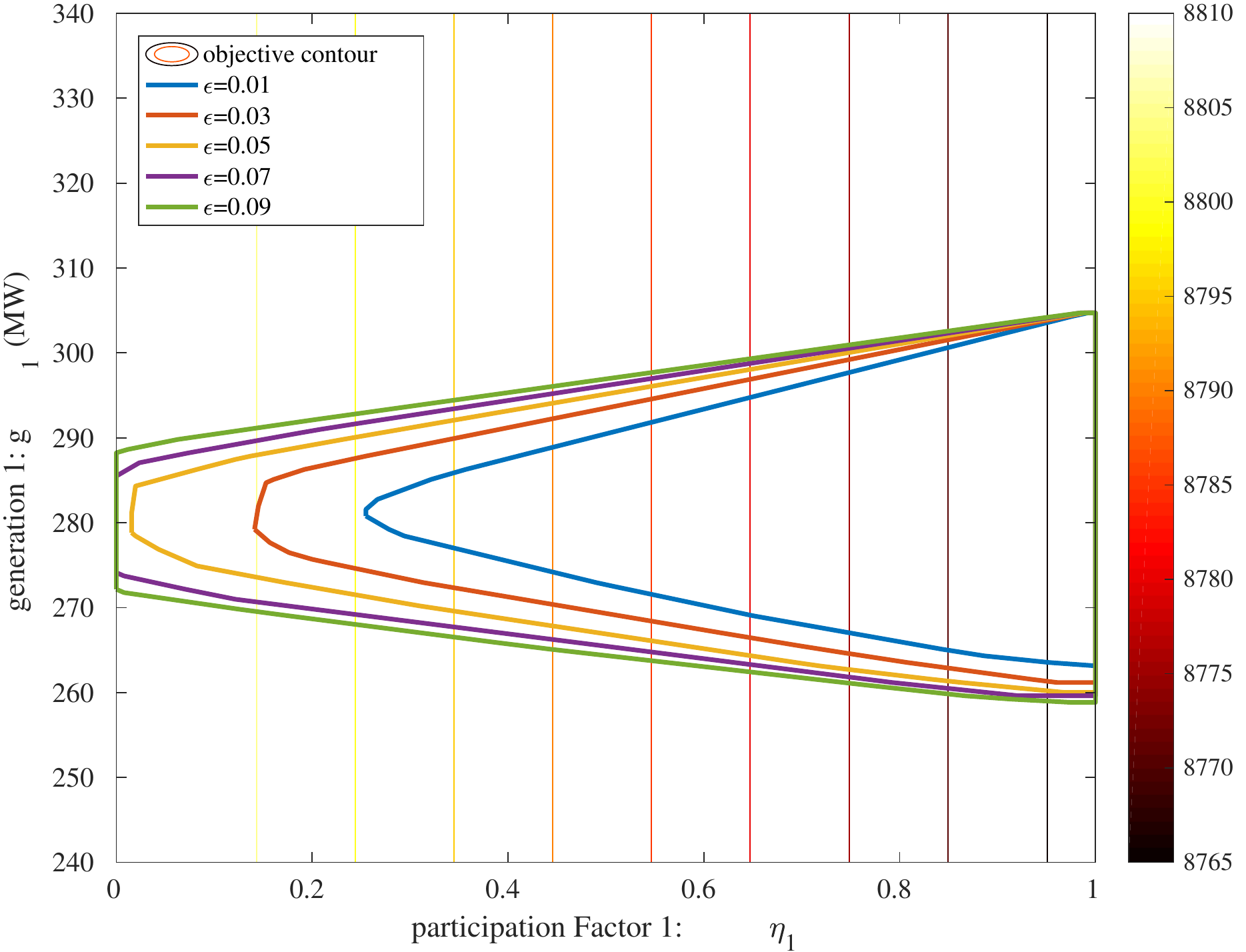}
	\caption{Feasible Region of cc-DCOPF for the 3-bus System (with contours of objective values)}
	\label{fig:ex_case3sc-feasible-region}
\end{figure}

Although there are four variables ($g_1,g_2,\eta_1,\eta_2$) in (cc-DCOPF) for the 3-bus system, only two of them (e.g. $g_1$ and $\eta_1$) are free variables because of constraints (\ref{form:cc-DCOPF-balance}) and (\ref{form:cc-DCOPF-sumone}) \footnote{$g_2 = \mathbf{1}^\intercal \hat{d} - g_1$ and $\eta_2 = 1-\eta_1$.}. Thanks to this, we could visualize the violation probability function $\mathbb{V}(g_1,\eta_1)$ in the 3D space. The function $\mathbb{V}(g_1,\eta_1)$ is evaluated over a grid of 65536 points in the $(g_1,\eta_1)$ space, the value of $\mathbb{V}(g_1,\eta_1)$ at every point is estimated using 1000 i.i.d realizations of uncertainties $\tilde{d}$. Figure \ref{fig:ex_case3sc-feasible-region-surf-view3} shows $\mathbb{V}(g_1,\eta_1)$. Based on the estimation of $\mathbb{V}(g_1,\eta_1)$, we could also visualize the feasible region $\mathcal{F}_\epsilon = \{x:\mathbb{V}(x) \le \epsilon\}$, which is shown in Figure \ref{fig:ex_case3sc-feasible-region}.

\subsection{Simulation Results} 
\label{sub:simulation_results}
We solve cc-DCOPF on the 3-bus system with eight different methods: (1) scenario approach with prior guarantees, (SA:prior, Corollary \ref{cor:prior_sample_complexity_fully_supported} in \citep{geng_data-driven_2019}); (2) scenario approach with posterior guarantees (SA:posterior, Theorem \ref{thm:posterior_guarantee} in \citep{geng_data-driven_2019}); (3) sample average approximation, where $N$ and $\varepsilon$ are chosen based on the sampling and discarding Theorem (SAA:s\&d, Theorem \ref{thm:sampling_and_discarding} in \citep{geng_data-driven_2019}); (4-7) Robust counterpart with different uncertainty sets specified in Theorem \ref{thm:rlo_safe_approx} in \citep{geng_data-driven_2019}: box (RC:box), ball (RC:ball), ball-box (RC:ball-box) and budget (RC:budget) uncertainty sets; (8) convex approximation with Markov bound (CA:markov, Theorem \ref{thm:convex_approx_safe_approx} and Proposition \ref{prop:convex_approx_cvar} in \citep{geng_data-driven_2019}). 

\begin{figure*}[htbp]
	\centering
	\begin{subfigure}[b]{0.49\textwidth}
    \includegraphics[width=\linewidth]{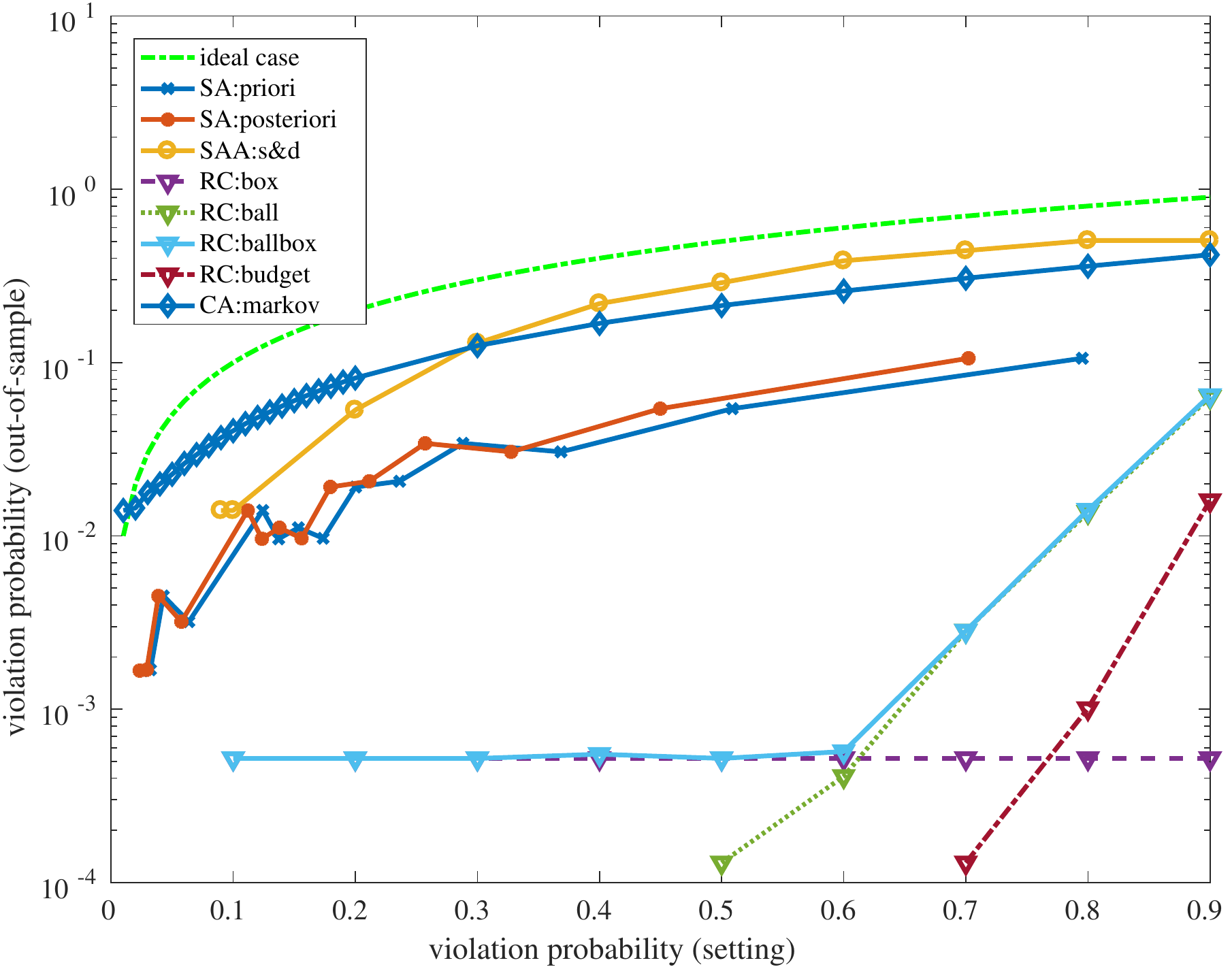}
    \caption{with logarithmic y-axis}
    \label{fig:ex_case3sc-all-methods-epsilon-normal}
    \end{subfigure}
    \begin{subfigure}[b]{0.49\textwidth}
    \includegraphics[width=\linewidth]{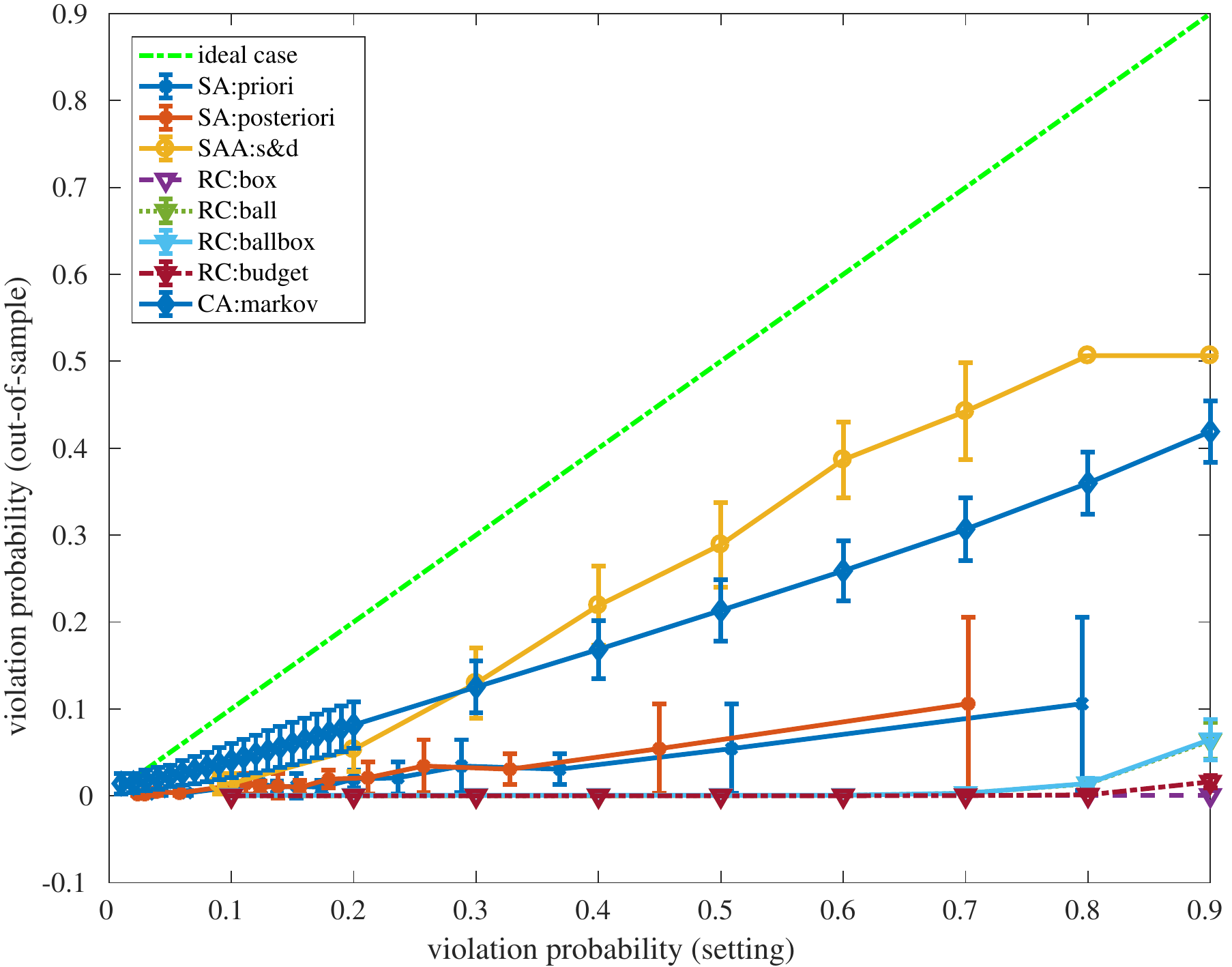}
    \caption{with error bars showing standard deviations}
    \label{fig:ex_case3sc-all-methods-epsilon-errorbar}
    \end{subfigure}
    \caption{Violation Probabilities (cc-DCOPF of the 3-bus System)}
    \label{fig:ex_case3sc-all-methods-epsilon-all}
\end{figure*}

\begin{figure*}[htbp]
	\centering
	\begin{subfigure}[b]{0.49\textwidth}
    \includegraphics[width=\linewidth]{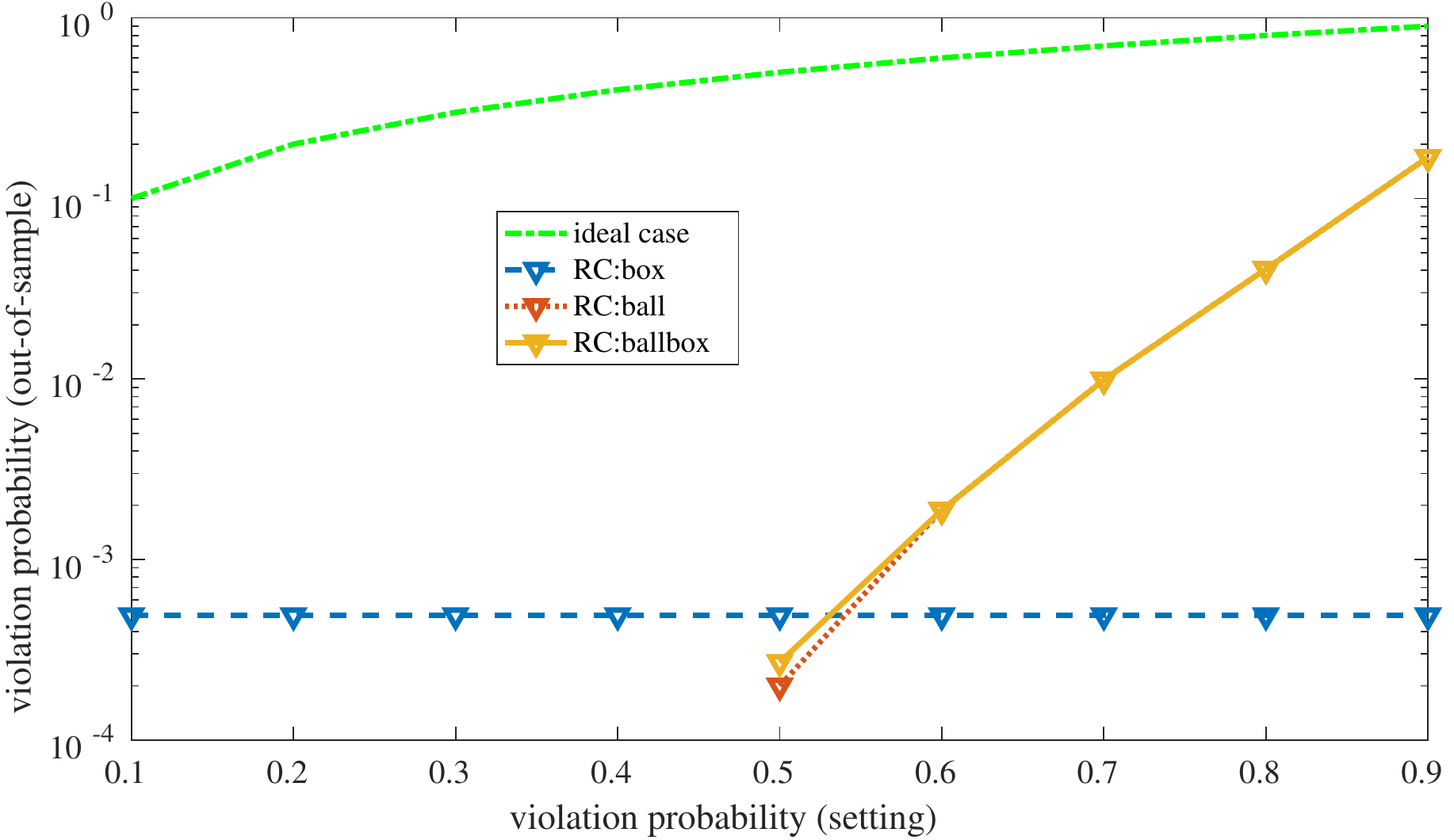}
    \caption{robust counterpart methods}
    \label{fig:ex_case24_ieee_rts-rc-methods-epsilon-gaussian}
    \end{subfigure}
    \begin{subfigure}[b]{0.49\textwidth}
    \includegraphics[width=\linewidth]{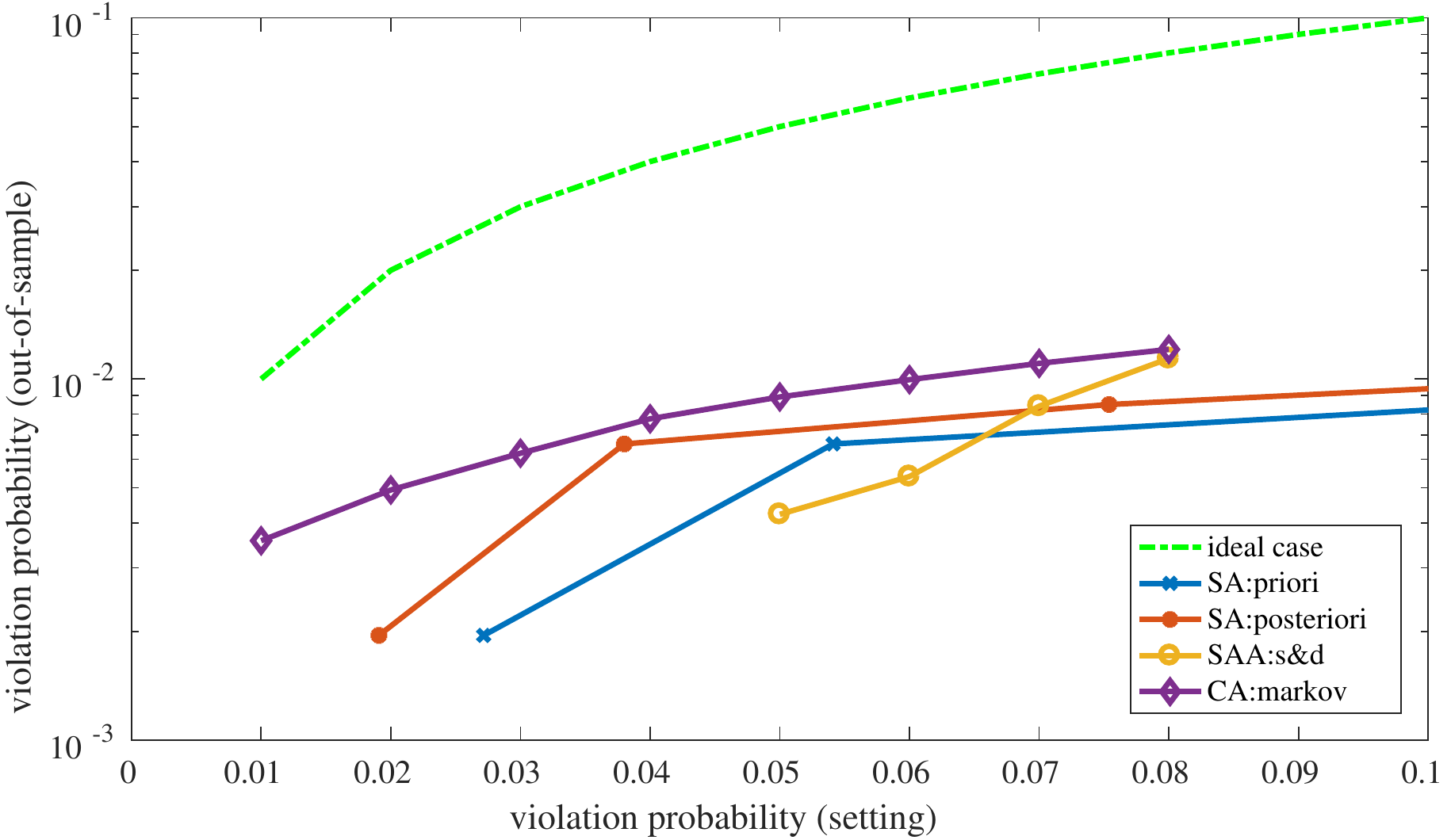}
    \caption{other methods}
    \label{fig:ex_case24_ieee_rts-all-methods-epsilon-gaussian}
    \end{subfigure}
    \caption{Violation Probabilities (cc-DCOPF of the 24-bus System, Gaussian Distributions)}
    \label{fig:ex_24_ieee_rts-all-methods-epsilon-all-gaussian}
\end{figure*}

We first examine the feasibility of the returned solutions from eight algorithms. Figure \ref{fig:ex_case3sc-all-methods-epsilon-all} and \ref{fig:ex_24_ieee_rts-all-methods-epsilon-all-gaussian} show the out-of-sample violation probabilities $\hat{\epsilon}$ versus desired $\epsilon$ in the setting. The green dashed lines in Figure \ref{fig:ex_case3sc-all-methods-epsilon-all} and \ref{fig:ex_24_ieee_rts-all-methods-epsilon-all-gaussian} denote the ideal case where $\hat{\epsilon} = \epsilon$. Any points above the green dashed line indicate infeasible solutions that $\mathbb{V}(x) > \epsilon$. Clearly all methods return feasible solutions (with high probability) to (CCO). From Figure \ref{fig:ex_case3sc-all-methods-epsilon-all}, sample average approximation and convex approximation are less conservative than other methods. However, it is worth noting that when $\epsilon$ is small (e.g. $10^{-2}$), the data-driven approximation of CVaR (Proposition \ref{prop:convex_approx_cvar}, \citep{geng_data-driven_2019}) does not necessarily give a safe approximation to (CCO) \citep{chen_cvar_2010}.
The robust counterpart methods are typically $10\sim100$ times more conservative than other methods, as illustrated in the comparison of Figure \ref{fig:ex_case24_ieee_rts-rc-methods-epsilon-gaussian} with Figure \ref{fig:ex_case24_ieee_rts-all-methods-epsilon-gaussian}. The conservativeness could be significantly reduced by better construction of uncertainty sets, e.g. \cite{chen_cvar_2010,bertsimas_data-driven_2018}. Among four different choices of uncertainty sets, the ball-box set is the least conservative one, which combines the advantages of the ball and box uncertainty sets.

\begin{figure}[tb]
	\centering
	\includegraphics[width=\linewidth]{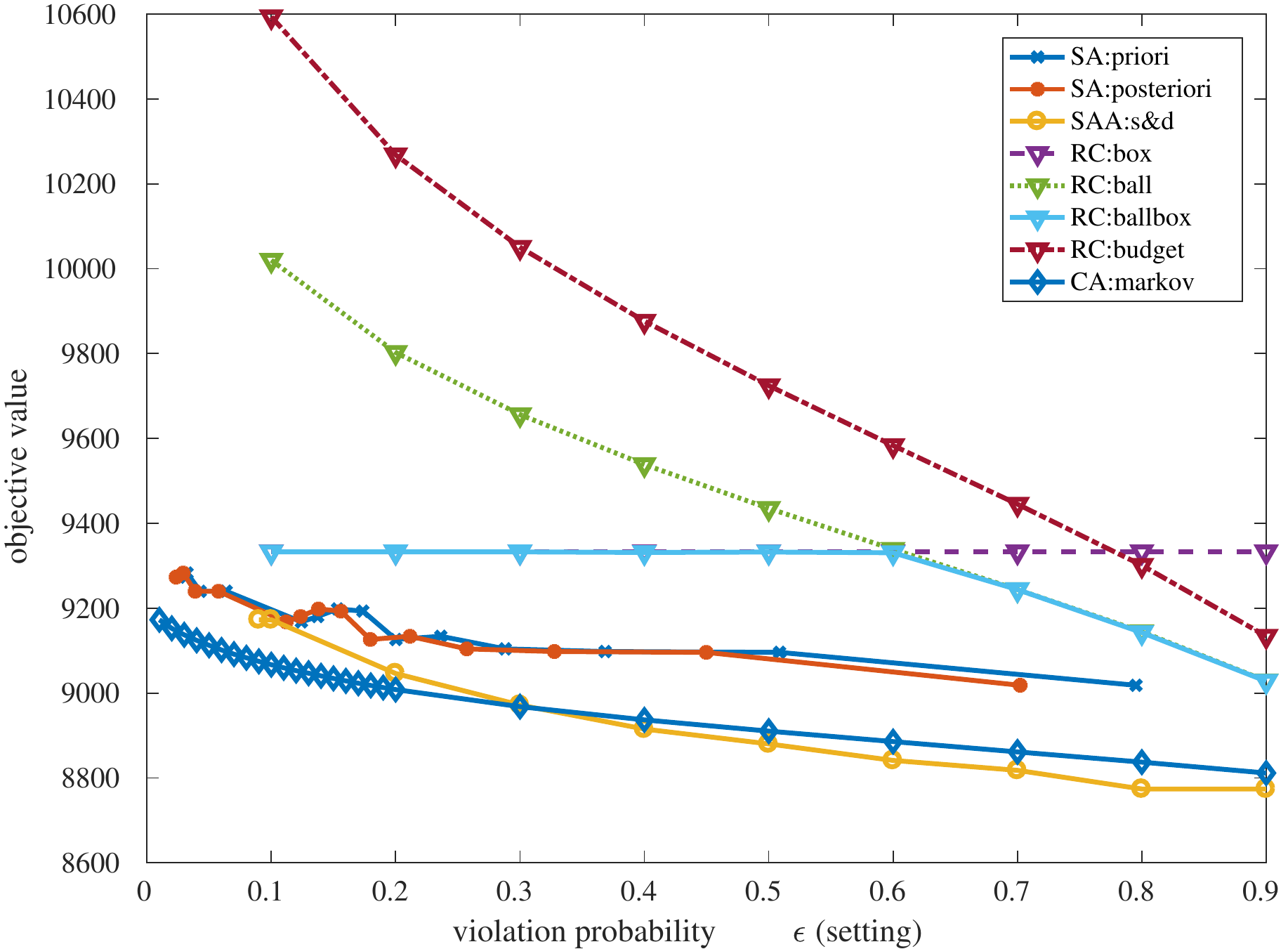}
	\caption{Objective Values (cc-DCOPF of the 3-bus System)}
	\label{fig:ex_case3sc-all-methods-objective}
\end{figure}

\begin{figure*}[htbp]
	\centering
	\begin{subfigure}[b]{0.49\textwidth}
    \includegraphics[width=\linewidth]{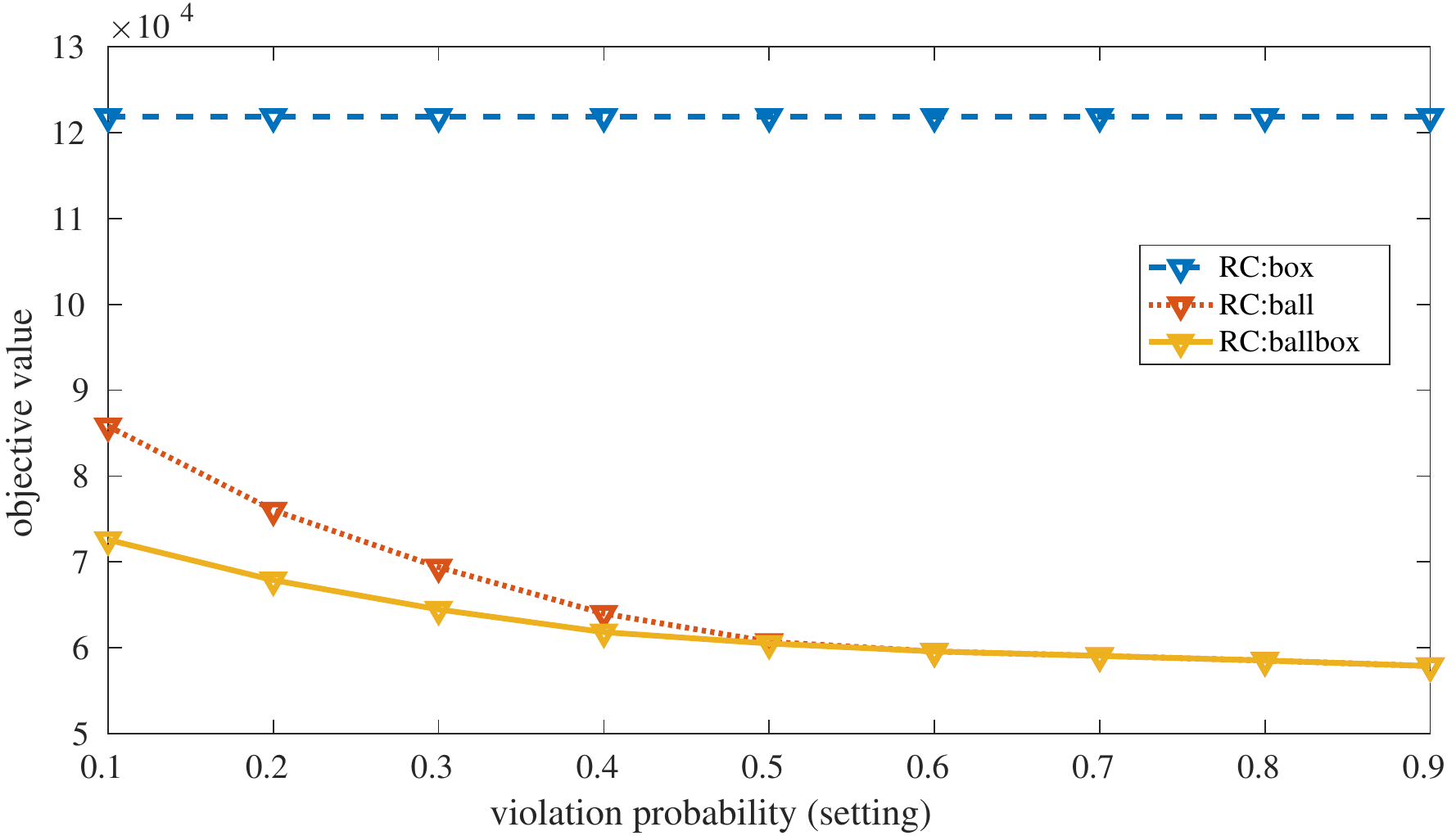}
    \caption{robust counterpart methods}
    \label{fig:ex_case24_ieee_rts-rc-methods-objective-gaussian}
    \end{subfigure}
    \begin{subfigure}[b]{0.49\textwidth}
    \includegraphics[width=\linewidth]{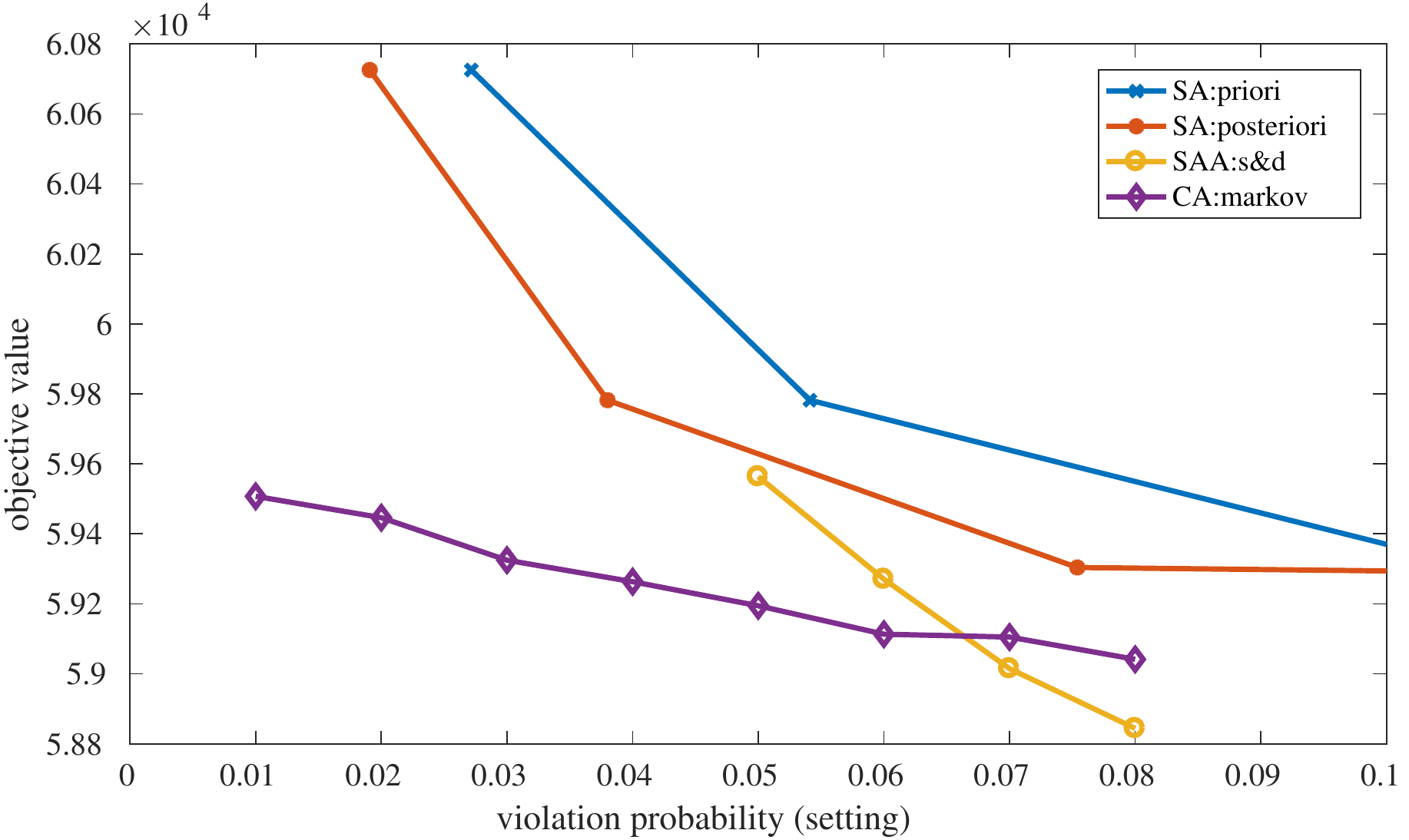}
    \caption{other methods}
    \label{fig:ex_case24_ieee_rts-all-methods-objective-gaussian}
    \end{subfigure}
    \caption{Objective Values (cc-DCOPF of the 24-bus System, Gaussian Distributions)}
    \label{fig:ex_case24_ieee_rts-all-methods-objective-all-gaussian}
\end{figure*}

\begin{figure}[htbp]
	\centering
    \includegraphics[width=\linewidth]{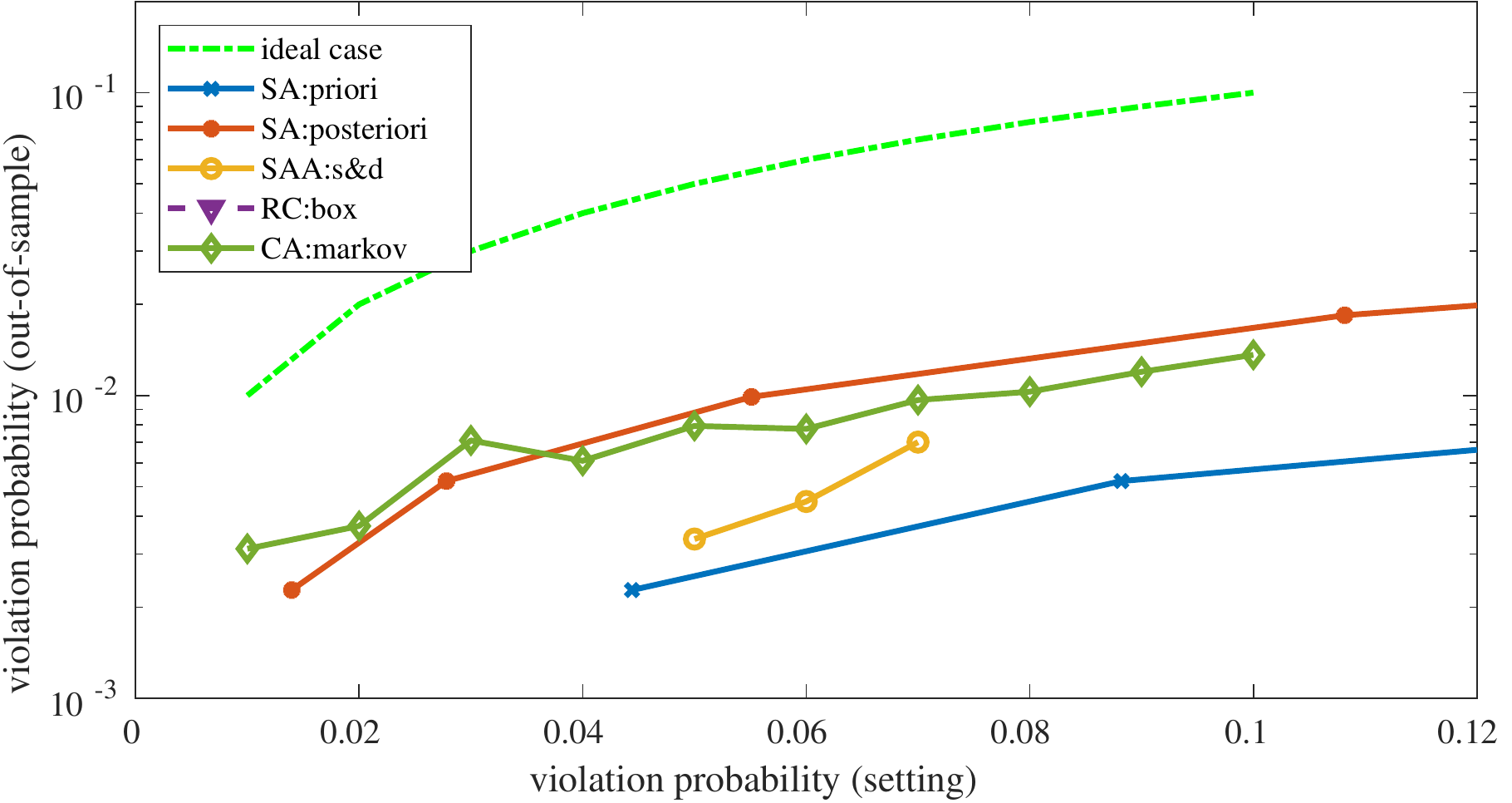}
    \caption{Violation Probabilities in Logarithmic Scale (cc-DCOPF of the 24-bus System, Beta Distributions)}
    \label{fig:ex_24_ieee_rts-all-methods-epsilon-all-beta}
\end{figure}

\begin{figure}[htbp]
	\centering
    \includegraphics[width=\linewidth]{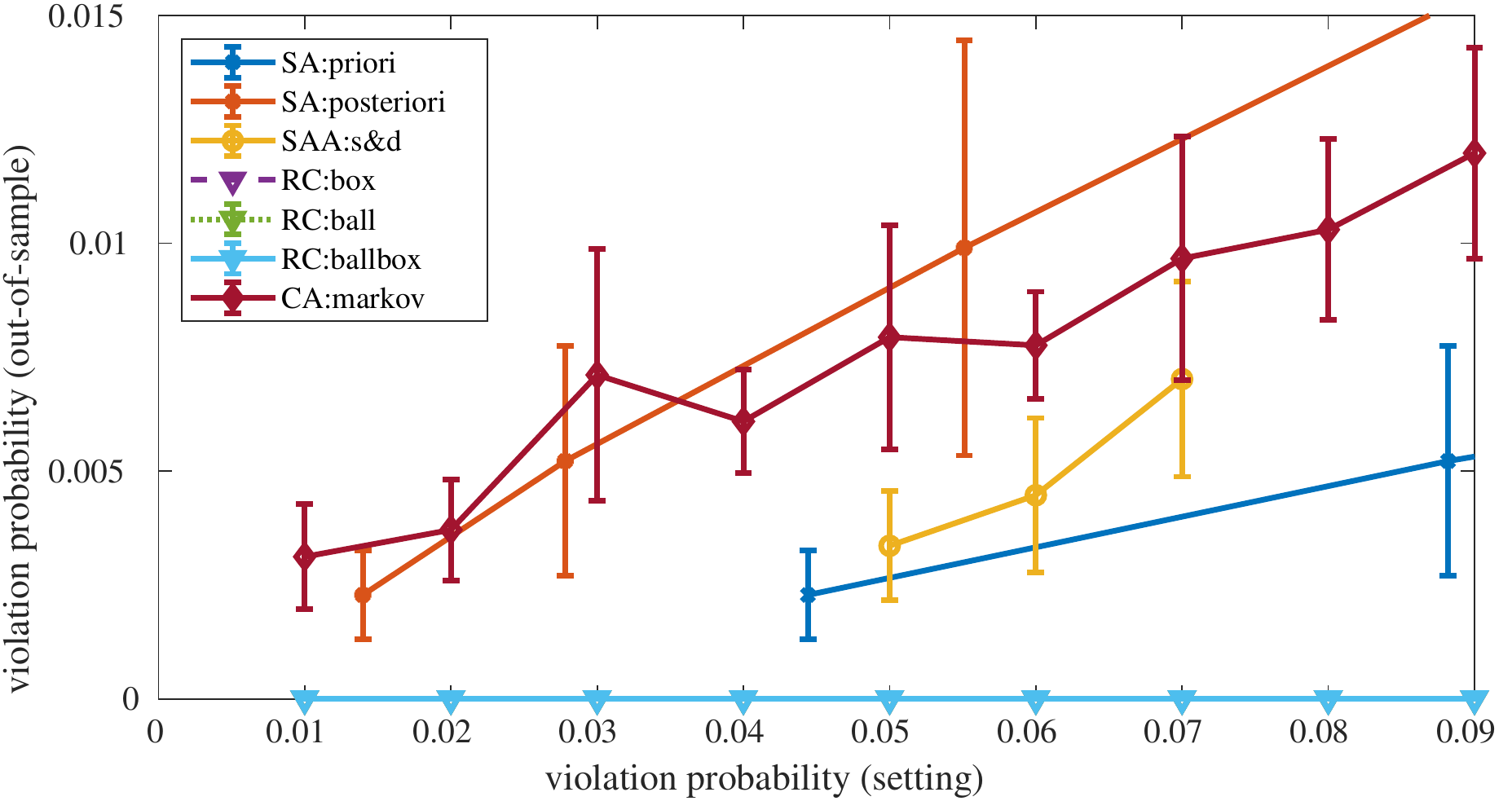}
    \caption{Violation Probabilities with error bars showing standard deviations (cc-DCOPF of the 24-bus System, Beta Distributions)}
    \label{fig:ex_24_ieee_rts-all-methods-epsilon-errorbar-all-beta}
\end{figure}

\begin{figure*}[htbp]
	\centering
	\begin{subfigure}[b]{0.49\textwidth}
    \includegraphics[width=\linewidth]{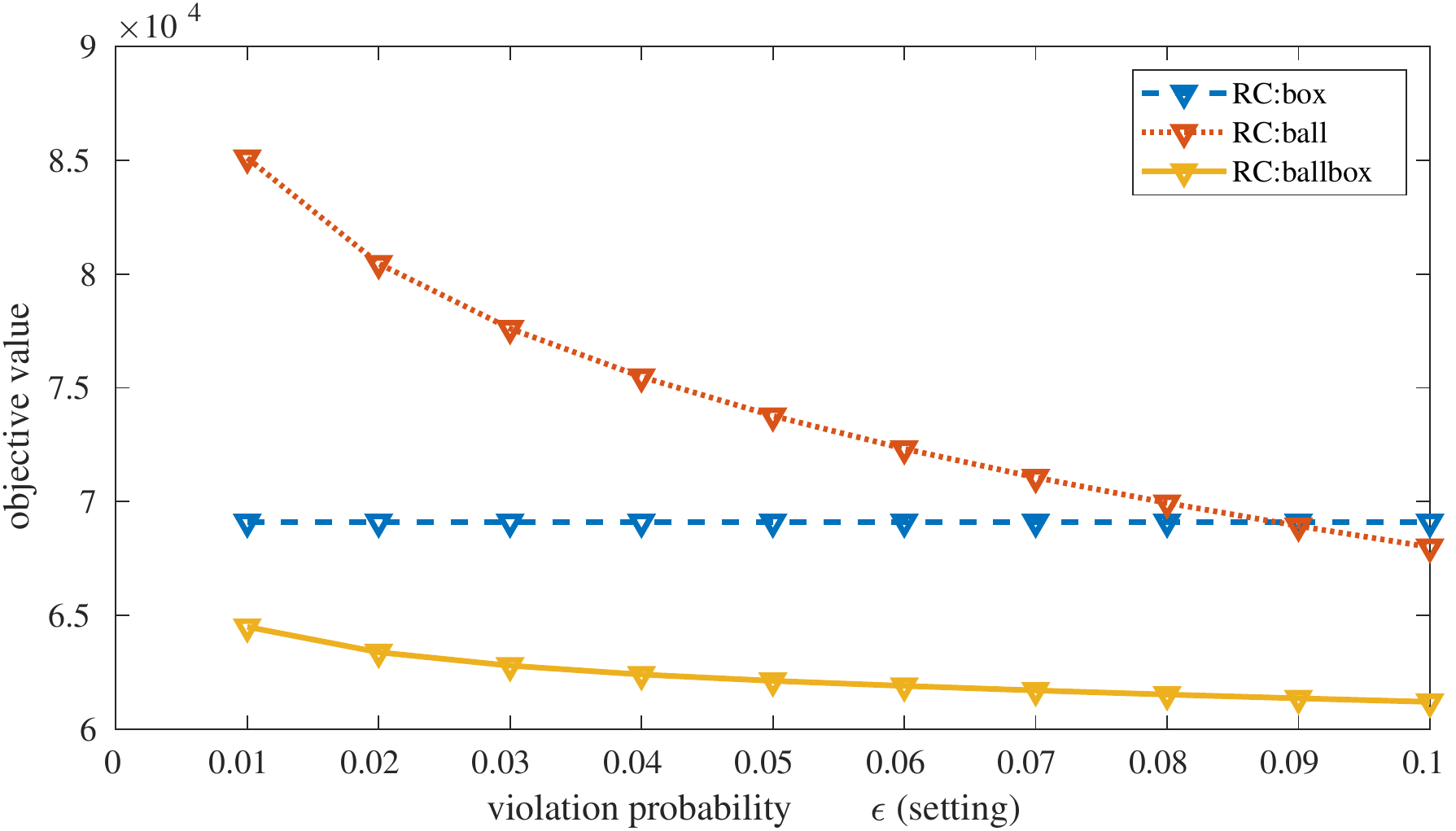}
    \caption{robust counterpart methods}
    \label{fig:ex_case24_ieee_rts-rc-methods-objective-beta}
    \end{subfigure}
    \begin{subfigure}[b]{0.49\textwidth}
    \includegraphics[width=\linewidth]{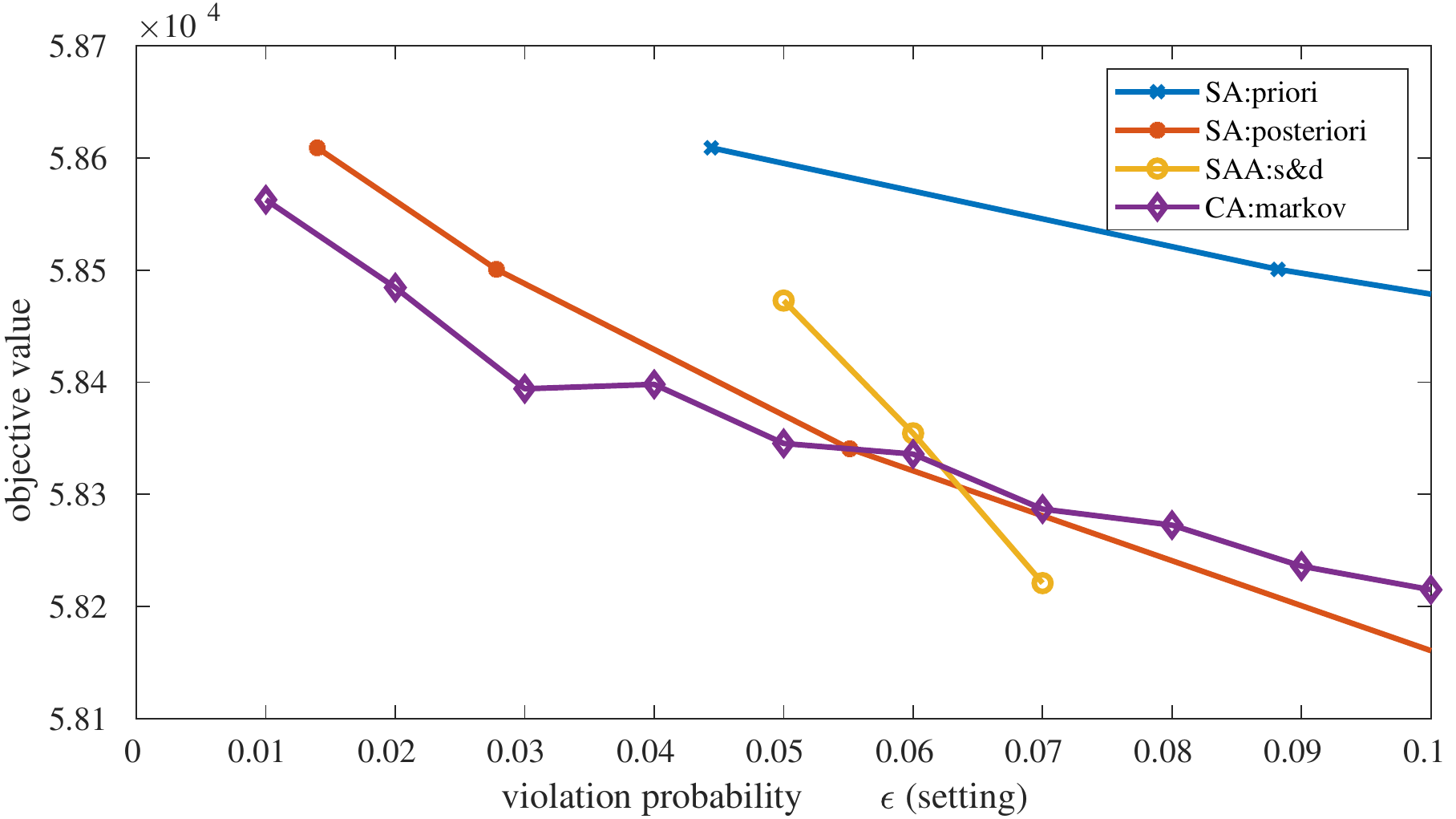}
    \caption{other methods}
    \label{fig:ex_case24_ieee_rts-all-methods-objective-beta}
    \end{subfigure}
    \caption{Objective Values (cc-DCOPF of the 24-bus System, Beta Distributions)}
    \label{fig:ex_case24_ieee_rts-all-methods-objective-all-beta}
\end{figure*}

Figure \ref{fig:ex_24_ieee_rts-all-methods-epsilon-all-gaussian}-\ref{fig:ex_case24_ieee_rts-all-methods-objective-all-gaussian} present the results of the 24-bus system with Gaussian distributions. Simulation results of the beta distribution are in Figure \ref{fig:ex_24_ieee_rts-all-methods-epsilon-all-beta}-\ref{fig:ex_case24_ieee_rts-all-methods-objective-all-beta}. Observations from Figure \ref{fig:ex_24_ieee_rts-all-methods-epsilon-all-beta}-\ref{fig:ex_case24_ieee_rts-all-methods-objective-all-beta} are similar with the case of Gaussian distributions. Every method behaves more conservative in the case of beta distributions than the case of Gaussian distributions. It is worth noting that the RO-based methods (RC:box, RC:ball, RC:ball-box in Figure \ref{fig:ex_24_ieee_rts-all-methods-epsilon-errorbar-all-beta}) are so conservative that lead to zero empirical violation probability $\hat{\epsilon}$. 





\section{Concluding Remarks} 
\label{sec:concluding_remarks}
This paper consists of two parts.
The first part presents a comprehensive review on the fundamental properties, key theoretical results, and three classes of algorithms for chance-constrained optimization. An open-source MATLAB toolbox ConvertChanceConstraint is developed to automate the process of translating chance constraints to compatible forms for mainstream optimization solvers. 
The second part of this paper presents three major applications of chance-constrained optimization in power systems. 
We also present a critical comparison of existing algorithms to solve chance-constrained programs on IEEE benchmark systems.

Many interesting directions are open for future research. More thorough and detailed comparisons of solutions to (CCO) on various problems with realistic datasets is needed. 
In terms of theoretical investigation, an analytical comparison of existing solutions to chance-constrained optimization is necessary to substantiate the fundamental insights obtained from numerical simulations. In terms of applications, many existing results can be improved by exploiting the structural properties of the problem to be solved. 
The application of chance-constrained optimization in electric energy systems could go beyond operational planning practices. For example, it would be worth investigating into the economic interpretation of market power issues through the lens of chance-constrained optimization.

\section*{References}
\bibliographystyle{elsarticle-harv}
\bibliography{references}

\end{document}